\documentclass[11pt,a4paper]{article}

\usepackage{graphicx}
\usepackage{amsmath}
\usepackage{amsthm}
\usepackage{amssymb}
\usepackage{amsfonts}
\usepackage{latexsym}

\newcommand{\beql}[1]{\begin{equation}\label{#1}}
\newcommand{\eeq}{\end{equation}}
\newcommand{\comment}[1]{}

\newcommand{\Norm}[1]{{\left\|{#1}\right\|}}

\newcommand{\PP}{{\mathcal P}}
\newcommand{\PK}{{\mathcal P}_n(K)}

\newcommand{\RR}{{\mathbb R}}
\newcommand{\CC}{{\mathbb C}}
\newcommand{\ZZ}{{\mathbb Z}}
\newcommand{\NN}{{\mathbb N}}

\newcommand{\diam}{{\rm diam\,}}

\newcounter{dfn}
\def\thedfn{\arabic{dfn}}

\newcounter{thm}
\newcounter{othm}
\def\theothm{\Alph{othm}}
\newenvironment{othm}{
  \em
  \vskip 0.10in
  \refstepcounter{othm}
  \noindent{\bf Theorem\ \theothm}
}{\vskip 0.10in}

\newcounter{mysec}
\def\themysec{\arabic{mysec}}
\newcommand{\mysection}[1]{
  \vskip 0.25in
  \refstepcounter{mysec}\centerline{\large\bf \S\themysec.\ {#1}}\par
  \nopagebreak
  \addcontentsline{toc}{section}{{\bf \themysec.}\ {#1}}
}
\newcounter{mysubsec}[mysec]
\newtheorem{theorem}{Theorem}

\newtheorem{lemma}{Lemma}

\newtheorem{definition}{Definition}

\theoremstyle{definition}

\newcounter{rem}
\newcounter{rep}
\reversemarginpar 

\begin{document}

\title{On a paper of Er\H od and Tur\'an-Markov inequalities for non-flat convex domains}
\author{Szil\' ard Gy. R\' ev\' esz}


\comment{\dedication{To the memory of J\'anos Er\H od, 1916-1945}}

\maketitle

{\small{\em To the memory of J\'anos Er\H od, 1916-1945}}

\vskip3cm

\begin{abstract}
For a convex domain $K\subset\CC$ the well-known general Markov
inequality asserting that a polynomial $p$ of degree $n$ must have
$\Norm{p'} \le c(K) n^2 \Norm{p}$ holds. On the other hand for
polynomials in general, $\Norm{p'}$ can be arbitrarily small as
compared to $\Norm{p}$.

The situation changes when we assume that the polynomials have all
their zeroes in the convex body $K$. This problem of lower bound
for Markov factors was first investigated by Tur\'an in 1939.
Tur\'an showed $\Norm{p'} \ge (n/2) \Norm{p}$ for the unit disk
$D$ and $\Norm{p'} \ge c \sqrt{n} \Norm{p}$ for the unit interval
$I:=[-1,1]$. Soon after that, J. Er\H od published a long article,
discussing various extensions of the results and methods of
Tur\'an.

For decades, Er\H od's paper was quoted only for the explicit
calculation of the exact constant of the interval case. However,
in recent years Levenberg and Poletsky, Erd\'elyi and also the
author investigated Tur\'an's problem for various sets --
basically, convex domains. In this context the much richer content
of Er\H od's work is to be realized again.

Thus, the aim of the paper is twofold. On the one hand we give an
account of the half-forgotten, old Hungarian article of Er\H od,
also commemorating its author. On the other hand we report on
recent developments with particular emphasis on development of one
of the key observations of Er\H od, namely, the role of the
\emph{curvature} of the boundary curve in the estimation of the
lower bound of Markov factors.
\end{abstract}

\vskip2cm



\let\oldfootnote\thefootnote
\def\thefootnote{}
\footnotetext{Supported in part in the framework of the
Hungarian-French Scientific and Technological Governmental
Cooperation, Project \# F-10/04 and the Hungarian-Spanish
Scientific and Technological Governmental Cooperation, Project \#
E-38/04.} \footnotetext{This work was completed during the
author's stay in Paris under his Marie Curie fellowship, contract
\# MEIF-CT-2005-022927.}
\let\thefootnote\oldfootnote

\newpage


\mysection{Introduction} \label{sec:introduction}


\noindent On the complex plane polynomials of degree $n$ admit a
Markov inequality $\Norm{p'}_K\le c_K n^2 \Norm{p}_K$ on all
convex, compact $K\subset \CC$. Here the norm $\Norm{\cdot}
:=\Norm{\cdot}_K$ denotes sup norm over values attained on $K$.

In 1939 Paul Tur\'an studied converse inequalities of the form
$\Norm{p'}_K\ge c_K n^A \Norm{p}_K$. Clearly such a converse can
hold only if further restrictions are imposed on the occurring
polynomials $p$. Tur\'an assumed that all zeroes of the
polynomials must belong to $K$. So denote the set of complex
(algebraic) polynomials of degree (exactly) $n$ as $\PP_n$, and
the subset with all the $n$ (complex) roots in some set
$K\subset\CC$ by $\PK$. The (normalized) quantity under our study
is thus the ``inverse Markov factor"
\begin{equation}\label{Mdef}
M_n(K):=\inf_{p\in \PK} M(p) \qquad \text{\rm with} \qquad
M:=M(p):=\frac{\Norm{p'}}{\Norm{p}}~.
\end{equation}

\begin{othm}{\bf[Tur\'an, \cite[p. 90]{Tur}].}\label{oth:Turandisk} If $p\in \PP_n(D)$, where
$D$ is the unit disk, then we have
\begin{equation}\label{Turandisk}
\Norm{p'}_D\ge \frac n2 \Norm{p}_D~.
\end{equation}
\end{othm}

\begin{othm}{\bf[Tur\'an, \cite[p. 91]{Tur}].}\label{oth:Turanint} If $p\in\PP_n(I)$,
where $I:=[-1,1]$, then we have
\begin{equation}\label{Turanint}
\Norm{p'}_I\ge \frac {\sqrt{n}}{6} \Norm{p}_I~.
\end{equation}
\end{othm}

Theorem \ref{oth:Turandisk} is best possible, as the example of
$p(z)=1+z^n$ shows. This also highlights the fact that, in
general, the order of the inverse Markov factor cannot be higher
than $n$. On the other hand, a number of positive results, started
with J. Er\H od's work, exhibited convex domains having order $n$
inverse Markov factors (like the disk). We come back to this after
a moment.

Regarding Theorem \ref{oth:Turanint}, Tur\'an pointed out by the
example of $(1-x^2)^{n}$ that the $\sqrt{n}$ order is sharp. The
slightly improved constant $1/(2e)$ can be found in \cite{LP}, but
the value of the constant is computed for all fixed $n$ precisely
in \cite{Er}. In fact, about two-third of the paper \cite{Er} is
occupied by the rather lengthy and difficult calculation of these
constants, which partly explains why later authors started to
consider this achievement the only content of the paper. Our aim
is to describe further ideas of Er\H od, as presented in
\cite{Er}, and to describe development of these ideas to date.

As mentioned above, Er\H od did not stop at calculation of
$M_n(I)$. He then considered ellipse domains, which form a
parametric family $E_b$ naturally connecting the two sets $I$ and
$D$. Note that for the same sets $E_b$ the best form of the
Bernstein-Markov inequality was already investigated by Sewell,
see \cite{Sew}.

\begin{othm}[Er\H od, \cite[p. 70]{Er}]\label{th:ellipse}
Let $0<b<1$ and let $E_b$ denote the ellipse domain with major
axes $[-1,1]$ and minor axes $[-ib,ib]$. Then
\begin{equation}\label{oldellipse}
\Norm{p'} \ge \frac {b}{2}  n \Norm{p}
\end{equation}
for all polynomials $p$ of degree $n$ and having all zeroes in
$E_b$.
\end{othm}

Er\H od himself provided two proofs, the first being a quite
elegant one using elementary complex functions, while the second
one fitting more in the frame of classical analytic geometry. In
2004 this theorem was rediscovered by J. Szabados, 
providing a testimony of the natural occurrence of the sets $E_b$
in this context\footnote{After learning about the overlap with
Er\H od's work, the result was not published.}.

In fact, the key to Theorem \ref{oth:Turandisk} was the following
observation, implicitly already in \cite{Tur} and \cite{Er} and
formulated explicitly in \cite{LP}.

\begin{lemma}[Tur\'an, Levenberg-Poletsky]\label{Tlemma} Assume
that $z\in\partial K$ and that there exists a disc $D_R$ of radius
$R$ so that $z\in\partial D_R$ and $K\subset D_R$. Then for all
$p\in\PK$ we have
\begin{equation}\label{Rdisc}
|p'(z)| \ge \frac n{2R} |p(z)|~.
\end{equation}
\end{lemma}

So Levenberg and Poletsky \cite{LP} found it worthwhile to
formally introduce the next definition.

\begin{definition}\label{circular} A compact set $K\subset \CC$ is
called $R$-{\em circular}, if for any point $z\in\partial K$ there
exists a disc $D_R$ of radius $R$ with $z\in\partial D_R$ and
$K\subset D_R$.
\end{definition}

With this they formulated various consequences. For our present
purposes let us chose the following form, c.f. \cite[Theorem
2.2]{LP}.

\begin{othm}[Er\H od; Levenberg-Poletsky]\label{th:circular} If $K$ is an
$R$-circular set and $p\in\PK$, then
\begin{equation}\label{circularineq}
\Norm{p'}\ge \frac n{2R}  \Norm{p}~.
\end{equation}
\end{othm}

Note that here it is not assumed that $K$ be convex; a circular
arc, or a union of disjoint circular arcs with proper points of
join, satisfy the criteria. However, other curves, like e.g. the
interval itself, do not admit such inequalities; as said above,
the order of magnitude can be as low as $\sqrt{n}$ in general.

Er\H od did not formulate the result that way; however, he was
clearly aware of that. This can be concluded from his various
argumentations, in particular for the next result.

\begin{othm}[Er\H od, \cite[p. 77]{Er}]\label{oth:smooth} If $K$ is a
$C^2$-smooth convex domain with the curvature of the boundary
curve staying above a fixed positive constant $\kappa>0$, and if
$p\in\PK$, then we have
\begin{equation}\label{circularestimate}
\Norm{p'}\ge c(K) n  \Norm{p}.
\end{equation}
\end{othm}

From Er\H od's argument one can not easily conclude that the
constant is $c(K)=\kappa/2$; on the other hand, his statement is
more general than that. Although the proof is slightly incomplete,
let us briefly describe the idea. We will return to this and
provide a somewhat different, complete proof giving also the value
$c(K)=\kappa/2$ of the constant later in \S \ref{sec:circular}.

\begin{proof} The norm of $p$ is attained at some point of the
boundary, so it suffices to prove that $|p'(z)|/|p(z)|\geq cn$ for
all $z\in \partial K$. But the usual form of the logarithmic
derivative and the information that all the $n$ zeroes
$z_1,\dots,z_n$ of $p$ are located in $K$ allows us to draw this
conclusion once we have for a fixed direction
$\varphi:=\varphi(z)$ the estimate
\begin{equation}\label{onerootgives}
\Re\left( e^{i\varphi}\frac{1}{z-z_k} \right)\geq c >0 \qquad
(k=1,\dots,n).
\end{equation}
Choosing $\varphi$ the (outer) normal direction of the convex
curve $\partial K$ at $z\in\partial K$, and taking into
consideration that $z_k$ are placed in $K\setminus \{z\}$
arbitrarily, we end up with the requirement that
\begin{equation}\label{onepointgives}
\Re\left( e^{i\varphi}\frac{1}{z-w}\right) =
\frac{\cos\alpha}{|z-w|} \geq c \qquad (w\in
K\setminus\{z\},\,\,\,\alpha:=\varphi-\arg(z-w))~.
\end{equation}
Now if $K$ is strictly convex, then for $z\ne w$ we do not have
$\cos\alpha=0$, a necessary condition for keeping the ratio off
zero. It remains to see if $|z-w|/\cos\alpha$ stays bounded when
$z\in\partial K$ and $w\in K\setminus\{z\}$, or, as is easy to
see, if only $w\in\partial K\setminus\{z\}$. Observe that
$F(z,w):=|z-w|/\cos\alpha$ is a two-variate function on $\partial
K^2$, which is not defined for the diagonal $w=z$, but under
certain conditions can be extended continuously. Namely, for given
$z$ the limit, when $w\to z$, is the well-known geometric quantity
$2\rho(z)$, where $\rho(z)$ is the radius of the osculating circle
(i.e., the reciprocal of the curvature $\kappa(z)$). (Note here a
gap in the argument for not taking into consideration also
$(z',w')\to (z,z)$, which can be removed by showing uniformity of
the limit.) Hence, for smooth $\partial K$ with strictly positive
curvature bounded away from 0, we can define
$F(z,z):=2/\kappa(z)=2\rho(z)$. This makes $F$ a continuous
function all over $\partial K^2$, hence it stays bounded, and we
are done.

To show the uniformity of the limit when $z',w'\to z$, let us fix
the arc length parametrization of $\gamma:=\partial K$, and assume
that $z=\gamma(s)$ and $w=\gamma(t)$: similarly $z'=\gamma(s')$
and $w'=\gamma(t')$. Now $w-z=\int_s^t\dot{\gamma}(u)du=
\int_s^t\int_s^u\ddot{\gamma}(v)dvdu + (t-s)\dot{\gamma}(s)$,
\begin{equation}\label{eqn:Fdef}
F(z,w)=\frac{\langle w-z ; \ddot{\gamma}(s)\rangle}
{|w-z|^2|\ddot{\gamma}(s)|}
\end{equation}
and so in view of $\langle \dot{\gamma}(s);\ddot{\gamma}(s)
\rangle =0$ we are led to
$$
F(z,w)=\frac{\int_s^t\int_s^u\langle
\ddot{\gamma}(v);\ddot{\gamma}(s)\rangle dvdu} {\langle
\int_s^t\int_s^u\ddot{\gamma}(v)dvdu + (t-s)\dot{\gamma}(s);
\int_s^t\int_s^u\ddot{\gamma}(v) dvdu +
(t-s)\dot{\gamma}(s)\rangle \left| \ddot{\gamma}(s)\right|}.
$$
In the numerator we can apply $\ddot{\gamma}\in C(\RR/L\ZZ)$
(where $L$ is the arc length of $\gamma$) and also $0<\kappa\leq
|\ddot{\gamma}| \leq \lambda <\infty$, say. Moreover, for $w$ in a
$\delta$-neighborhood of $z$ (in arc length distance), we even
have $|\ddot{\gamma}(v)-\ddot{\gamma}(s)|<\epsilon$ and hence the
numerator can be reformulated as
$$
\frac 12 (t-s)^2 |\ddot{\gamma}(s)|^2 + \int_s^t\int_s^u
\delta(s,v) dvdu = \frac 12 (t-s)^2 (|\ddot{\gamma}(s)|^2 +
\eta(s,t)) \qquad (|\delta(s,v)|, |\eta(s,t)| <\epsilon).
$$
On the other hand already $ (0<) |\ddot{\gamma}|\leq \lambda $
suffices (using also the fact that $|\dot{\gamma}(s)|=1$ in arc
length parametrization) to get for the denominator that it is
$$
(t-s)^2 \left(1+ \theta(s,t) \right)|\ddot{\gamma}(s)| \qquad
\left(|\theta(s,t)| < \delta \lambda + \frac 14 \delta^2 \lambda^2
\right).
$$
Summing up, we are led to
$$
F(z,w)=\frac{\frac 12 (t-s)^2 (|\ddot{\gamma}(s)|^2 (1+o(1))}
{(t-s)^2 \left( 1+o(1)\right)
|\ddot{\gamma}(s)|}=\frac{1+o(1)}{2}|\ddot{\gamma}(s)|\qquad (w\to
z).
$$
Note that by the uniform continuity of $\ddot{\gamma}$, the $o(1)$
is uniform (in arc length parameter); therefore for the pair of
points $(z',w')$ we similarly obtain
$$
F(z',w')=(1+o(1))\frac 12|\ddot{\gamma}(s')|=(1+o(1))\frac
12|\ddot{\gamma}(s)|
$$
if $(z',w')\to z$, as needed.
\end{proof}

From this argument it can be seen that whenever we have the
property \eqref{onepointgives} for all given boundary points $z\in
\partial K$, then we also conclude the statement. This
explains why Er\H od could allow even vertices, relaxing the
conditions of the above statement to hold only piecewise on smooth
Jordan arcs, joining at vertices. However, to have a fixed bound,
either the number of vertices has to be bounded, or some
additional condition must be imposed on them. Er\H od did not
elaborate further on this direction.

Convex domains (or sets) {\em not} satisfying the $R$-circularity
criteria with any fixed positive value of $R$ are termed to be
{\em flat}. Clearly, the interval is flat, like any polygon or any
convex domain which is not strictly convex. From this definition
it is not easy to tell if a domain is flat, or if it is circular,
and if so, then with what (best) radius $R$. We will deal with the
issue in this work, aiming at finding a large class of domains
having $cn$ order of the inverse Markov factor with some
information on the arising constant as well.

On the other hand a lower estimate of the inverse Markov factor of
the same order as for the interval was obtained in full generality
in 2002, see \cite[Theorem 3.2]{LP}.

\begin{othm}[Levenberg-Poletsky]\label{th:generalroot} If
$K\subset \CC$ is a compact, convex set, $d:=\diam{K}$ is the
diameter of $K$ and $p\in \PK$, then we have
\begin{equation}\label{genrootineq}
\Norm{p'}\ge \frac {\sqrt{n}}{20\,\diam(K)}  \Norm{p}~.
\end{equation}
\end{othm}

Clearly, we can have no better order, for the case of the interval
the $\sqrt{n}$ order is sharp. Nevertheless, already Er\H od
\cite[p. 74]{Er} addressed the question: ``For what kind of
domains does the method of Tur\'an apply?" Clearly, by ``applies"
he meant that it provides $cn$ order of oscillation for the
derivative.

The most general domains with $M(K)\gg n$, found by Er\H od, were
described on p. 77 of \cite{Er}. Although the description is a bit
vague, and the proof shows slightly less, we can safely claim that
he has proved the following result.

\begin{othm}[Er\H od\label{th:transfquarter} Let $K$ be any convex domain bounded by finitely many
Jordan arcs, joining at vertices with angles $<\pi$, with all the
arcs being $C^2$-smooth and being either straight lines of length
$\ell<\Delta(K)/4$, where $\Delta(K)$ stands for the transfinite
diameter of $K$, or having positive curvature bounded away from
$0$ by a fixed constant. Then there is a constant $c(K)$, such
that $M_n(K)\geq c(K) n$ for all $n\in\NN$.
\end{othm}

To deal with the flat case of straight line boundary arcs, Er\H od
involved another approach, cf. \cite[p. 76]{Er}, appearing later
to be essential for obtaining a general answer. Namely, he quoted
Faber \cite{Faber} for the following fundamental result going back
to Chebyshev.

\begin{lemma}[{\bf Chebyshev}]\label{l:capacity} Let $J=[u,v]$ be any interval on
the complex plane with $u\ne v$ and let $J \subset R \subset \CC$
be any set containing $J$. Then for all $k\in\NN$ we have
\begin{equation}\label{capacity}
\min_{w_1,\dots,w_k\in R} \max_{z\in J} \left| \prod_{j=1}^k
(z-w_j) \right| \ge 2 \left(\frac{|J|}{4}\right)^k~.
\end{equation}
\end{lemma}

The relevance of Chebyshev's Lemma is that it provides a
quantitative way to handle contribution of zero factors at some
properly selected set $J$. One uses this for comparison: if
$|p(\zeta)|$ is maximal at $\zeta\in\partial K$, then the maximum
on some $J$ can not be larger. Roughly speaking, combining this
with geometry we arrive at an effective estimate of the
contribution, hence even on the location of the zeroes. For more
in this direction see \cite{Er,R}.

In his recent work \cite{E2}, Erd\'elyi considered various special
domains. Apart from further results for polynomials of some
special form (e.g. even or real polynomials), he obtained the
following.

\begin{othm}[Erd\'elyi]\label{th:square}
Let $Q$ denote the square domain with diagonal $[-1,1]$. Then for
all polynomials $p\in\PP_n(Q)$ we have
\begin{equation}\label{oldsquare}
\Norm{p'} \ge C_0  n \Norm{p}
\end{equation}
with a certain absolute constant $C_0$.
\end{othm}

Note that the regular $n$-gon $K_n$ is already covered by Er\H
od's Theorem \ref{th:transfquarter} if $n\geq 26$, but not the
square $Q$, since the side length $h$ is larger than the quarter
of the transfinite diameter $\Delta$: actually, $\Delta(Q)\approx
0.59017\dots h$, while
$$
\Delta(K_n)=
\frac{\Gamma(1/n)}{\sqrt{\pi}2^{1+2/n}\Gamma(1/2+1/n)} h > 4h
\quad{\hbox{iff}}\qquad n\geq 26,
$$
{{see \cite[p. 135]{Rans}}}. Erd\'elyi's proof is similar to Er\H
od's argument\footnote{Erd\'elyi was apparently not aware of the
full content of \cite{Er} when presenting his rather similar
argument.}: sacrificing generality gives the possibility for a
better calculation for the particular choice of $Q$.

Returning to the question of the order in general, let us recall
that the term {\em convex domain} stands for a compact, convex
subset of $\CC$ \emph{having nonempty interior}. Clearly, assuming
boundedness is natural, since all polynomials of positive degree
have $\Norm{p}_K=\infty$ when the set $K$ is unbounded. Also, all
convex sets with nonempty interior are {\em fat}, meaning that
${\rm cl}(K)={\rm cl}({\rm int} K)$. Hence taking the closure does
not change the $\sup $ norm of polynomials under study. The only
convex, compact sets, falling out by our restrictions, are the
intervals, for what Tur\'an has already shown that his $c\sqrt{n}$
lower estimate is of the right order. Interestingly, it turned out
that among all convex compacta only intervals can have an inverse
Markov constant of such a small order.

\begin{othm}[Hal\'asz and R\'ev\'esz, \cite{R}]\label{th:cn} Let
$K\subset\CC$ be any convex domain having minimal width $w(K)$ and
diameter $d(K)$. Then for all $p\in\PK$ we have
\begin{equation}\label{cnresult}
\frac{\Norm{p'}}{\Norm{p}} \ge C(K) n \qquad \text{\rm with}
\qquad C(K)=0.0003 \frac{w(K)}{d^2(K)}~.
\end{equation}
\end{othm}

In the proof of this result in \cite{R}, due to generality, the
precision of constants could not be ascertained e.g. for the
special ellipse domains considered in \cite{Er}. Thus it seems
that the general results are not capable to fully cover e.g.
Theorem \ref{th:ellipse}.

Our aim here is to show that even that is possible for a quite
general class of convex domains with order $n$ inverse Markov
factors and a different estimate of the arising constants. This
will be achieved working more in the direction of Er\H od's first
observation, i.e. utilizing information on curvature. Since these
results need some technical explanations, in particular for the
geometric terms we use, formulation of these will be postponed
until \S \ref{sec:circular}. Before that, the next section is
dedicated to the life and work of J\'anos Er\H od, and in \S
\ref{sec:geometry} we start with describing the underlying
geometry.


\mysection{A few words about the life and work of J\'anos Er\H od}
\label{sec:erod}

With this paper we would like to call the attention of the
approximation theory community to the rich content of the original
paper \cite{Er}. It is necessary since out of the dozen or so
references in the literature to \cite{Er}, none of these works
mention -- and, actually, very few people are aware of the fact --
that Er\H od's work covered a lot more than the mere calculation
of $M_n(I)$. The paper was written and published in Hungarian,
back in the eve of World War II, and in spite of the fact that
both the Mathematical Reviews and the Zentralblatt reviews mention
the general features of the paper, that aspect seems to be
forgotten. A particular aim of our paper is to commemorate J\'anos
Er\H od, the person, too.

J\'anos Er\H od was born to the Ehrlich family in Gy\"ongy\"os, a
city some 80 km East-East-North of Budapest, on 30 November 1916,
during World War I. The Jewish family had three children: J\'anos
was born second, between his two sisters Carmen and M\'arta.
Sometimes in the 1920's the family converted to the protestant
church; on this occasion the family name was changed to the
Hungarian name "Er\H od", although the parents kept the name
"Ehrlich". J\'anos learnt very well and graduated with an
excellent grade; furthermore, he was a successful problem solver
of the legendary "K\"oMaL, K\"oz\'episkolai Matematikai \'es
Fizikai Lapok" ("Secondary School Mathematics and Physics
Journal"). Therefore, he continued studies in mathematics and
physics at the Budapest University of Sciences. He was only 23
when he received his PhD in mathematics in 1939: his thesis is
just the reprint of the only paper \cite{Er} he wrote. Although
the topic is a continuation of the work of Paul Tur\'an, it can
not be seen from the dry quotations how close personal contacts
they might have had. Nevertheless, Tur\'an and Er\H od mutually
refer to each other in \cite{Tur} and \cite{Er}, so at least they
knew about each other.

Because of the Jewish laws already in effect, he could not hope
for a university employment. However, he registered to the
Reformed Church Theology College in P\'apa, another city about 120
km to the West from Budapest. Also there he graduated with
excellent grade after completing the four year curriculum in the
three years 1939-1942. He passed his first and second clergyman
exams in 1943 and 1944, again with excellence. Becoming a reformed
church priest, he could serve his church at various locations
including the vicinity of P\'apa and Gy\H or. For a while he
became the director of the church's orphan boys' house in
Kom\'arom, some 80 km's West-North-West of Budapest (now belonging
to Slovakia).

In 1944 his parents and his younger sister M\'arta were deported.
They were taken to Auschwitz - none of them returned. In February,
1945 J\'anos decided to return to P\'apa to his fiance, Jol\'an
Nemes. He could stay unnoticed only for a very short time. He was
arrested together with Jol\'anka. The cause formally was not that
he was a Jew, but some (rather unrealistic) accusations of
treachery by establishing a radio contact with the advancing
Soviet troops. The young couple was interrogated in the military
base of P\'apa. Dezs\H o Tr\'ocs\'anyi, J\'anos Er\H od's theology
professor of the College, protested against the brutal torture of
J\'anos, but to no avail. The young couple was killed, very likely
in the barracks. Their remnants were not found. Neither their
grave, nor the exact date of their death is known.

However, the mathematical achievements of J\'anos were not lost,
even if somewhat forgotten. It is in order to commemorate also its
martyr author, when reflecting back to the rich content of this
pioneering work.


\mysection{Some geometrical notions} \label{sec:geometry}

Recall that the term {\em convex domain} stands for a compact,
convex subset of $\CC\cong\RR^2$ having nonempty interior. For a
convex domain $K$ any interior point $z$ defines a parametrization
$\gamma(\varphi)$ of the boundary $\partial K$, taking the unique
point \hbox{$\{z+te^{i\varphi}:\,t\in (0,\infty)\}\cap \partial
K$} for the definition of $\gamma(\varphi)$. This defines the
closed Jordan curve $\Gamma=\partial K$ and its parametrization
$\gamma : [0,2\pi] \to \CC$. By convexity, at any boundary point
$\zeta=\gamma(\theta)\in \partial K$, the chords to boundary
points in some small vicinity of $\zeta$ with parameter $<\theta$
or with $>\theta$ have arguments below and above the argument of
the direction of any tangential (supporting) line at $\zeta$. Thus
the tangent direction or argument function $\alpha_{-}(\theta)$
can be defined as e.g. the supremum (or $\limsup$) of arguments of
chords from the left; similarly, $\alpha_{+}(\theta):=\inf \{\arg
(z-\zeta)~:~ z=\gamma(\varphi),~ \varphi>\theta \}$, and any line
$\zeta+e^{i\beta}\RR$ with $\alpha_{-}(\theta)\le \beta \le
\alpha_{+}(\theta)$ is a supporting line to $K$ at
$\zeta=\gamma(\theta)\in\partial K$. In particular the curve
$\gamma$ is differentiable at $\zeta=\gamma(\theta)$ if and only
if $\alpha_{-}(\theta)=\alpha_{+}(\theta)$; in this case the
tangent of $\gamma$ at $\zeta$ is $\zeta+e^{i\beta}\RR$ with the
unique value of $\alpha=\alpha_{-}(\theta)=\alpha_{+}(\theta)$. It
is clear that interpreting $\alpha_{\pm}$ as functions on the
boundary points $\zeta\in\partial K$, we obtain a
parametrization-independent function. In other words, we are
allowed to change parameterizations to arc length, say, when in
case of $|\Gamma|=a$ with $\Gamma=\partial K$ the functions
$\alpha_{\pm}$ map from $[0,a]$ to $[0,2\pi]$.

Observe that $\alpha_{\pm}$ are nondecreasing functions with total
variation ${\rm Var}\,[\alpha_{\pm}] = 2\pi$, and that they have a
common value precisely at continuity points, which occur exactly
at points where the supporting line to $K$ is unique. At points of
discontinuity $\alpha_{\pm}$ is the left-, resp. right continuous
extension of the same function. For convenience, and for better
matching with \cite{BS}, we may even define the function
$\alpha:=(\alpha_{+}+\alpha_{-})/2$ all over the parameter
interval.

For obvious geometric reasons we call the jump function
$\beta:=\alpha_{+}-\alpha_{-}$ the {\em supplementary angle}
function. In fact, $\beta$ and the usual Lebesgue decomposition of
the nondecreasing function $\alpha_{+}$ to
$\alpha_{+}=\sigma+\alpha_{*}+\alpha_{0}$, consisting of the pure
jump function $\sigma$, the nondecreasing singular component
$\alpha_{*}$, and the absolute continuous part $\alpha_0$, are
closely related. By monotonicity there are at most countable many
points where $\beta(x)>0$, and in view of bounded variation we
even have $\sum_x \beta(x) \le 2\pi$, hence the definition
$\mu:=\sum_x \beta(x)\delta_x$ defines a bounded, positive Borel
measure. Now it is clear that $\sigma(x)=\mu([0,x])$, while
$\alpha_{*}'=0$ a.e., and $\alpha_0$ is absolutely continuous. In
particular, $\alpha$ or $\alpha_{+}$ is differentiable at $x$
exactly when $\beta(x)=0$ and $x$ is not in the exceptional set of
non-differentiable points with respect to $\alpha_{*}$. That is,
we have differentiability almost everywhere, and
\begin{align}\label{aedifferentiability}
\int_x^y \alpha'(t) dt =& \alpha_0(y)-\alpha_0(x)\notag \\=&
[\alpha_{+}(y)-\sigma(y)-\alpha_{*}(y))]-[\alpha_{+}(x)-\sigma(x)-\alpha_{*}(x)]
\\ \le & \alpha_{-}(y)-\alpha_{+}(x)~.\notag
\end{align}
It follows that we have the criteria
\begin{equation}\label{differentiallarge}
\alpha'(t)\ge \lambda \qquad \text{a.e}. \quad t\in [0,a]
\end{equation}
if and only if
\begin{equation}\label{fixchange}
\alpha_{\pm}(y)-\alpha_{\pm}(x) \ge \lambda (y-x) \qquad \forall
x,y \in [0,a]~.
\end{equation}
Here we reserved to the arc length parametrization. Recall that
one of the most important geometric quantities, curvature, is just
$\kappa(s):=\alpha'(s)$, whenever parametrization is by arc length
$s$.

Thus we can rewrite \eqref{differentiallarge} as
\begin{equation}\label{curvaturesmall}
\kappa(t) \ge \lambda \qquad \text{a.e}. \quad t\in [0,a]~,
\end{equation}
or, with radius of curvature $\rho(t):=1/\kappa(t)$ introduced,
\begin{equation}\label{curvradlarge}
\rho(t) \le \frac{1}{\lambda} \qquad \text{a.e}.\quad t\in [0,a]~.
\end{equation}
Again, $\rho$ is a parametrization-invariant quantity (describing
the radius of the osculating circle). Actually, it is easy to
translate all these conditions to arbitrary parametrization of the
tangent angle function $\alpha$. Since also curvature and
curvature radius are parametrization-invariant quantities, all the
above hold for any parametrization.

Moreover, with a general parametrization let
$|\Gamma(\eta,\zeta)|$ stand for the arc length of the rectifiable
Jordan arc $\Gamma(\eta,\zeta)$ of the curve $\Gamma$ between the
two points $\zeta, \eta \in \Gamma=\partial K$. We can then say
that the curve satisfies a Lipschitz-type increase or {\em
subdifferential condition} whenever
\begin{equation}\label{subdiffcond}
|\alpha_{\pm}(\eta)-\alpha_{\pm}(\zeta)| \ge \lambda
|\Gamma(\eta,\zeta)| \qquad (\forall \zeta, \eta \in \Gamma)~.
\end{equation}
Clearly, the above considerations show that all the above are
equivalent.

In the paper we use the notation $\alpha$ (and also
$\alpha_{\pm}$) for the tangent angle, $\kappa$ for the curvature,
and $\rho$ for the curvature radius. These notations we will use
basically in function of the arc length parametrization $s$, but
with a slight abuse of notation also $\alpha_{-}(\varphi)$,
$\kappa(\zeta)$ etc. may occur with the obvious meaning.

\mysection{Results for non-flat domains} \label{sec:circular}

The above Theorem \ref{th:ellipse} was formulated with very
precise constants. In particular, it gives a good description of
the "inverse Markov factor"
$$
M(E_b):= \inf_{p\in\PP_n(E_b)} M(p),
$$
when $n$ is fixed and $b\to 0$. In this section we aim at a
precise generalization of Theorem \ref{th:ellipse} using
appropriate geometric notions. Our argument stems out of the
notion of "circular sets", used in \cite{LP} and going back to
Tur\'an's work. This approach
\comment{yields the value $M_n(E_b)\geq (b/2) n$ for all ellipse
domains $E_{b}$ with $0<b<1$. Thus we}
can indeed cover the full content of Theorem \ref{th:ellipse}.
Moreover, the geometric observation and criteria we present will
cover a good deal of different, not necessarily smooth domains.
First let us have a recourse to Theorem \ref{oth:smooth}.

\begin{theorem}\label{th:twicesmooth} Let $K\subset\CC$ be any convex
domain with $C^2$-smooth boundary curve $\partial K=\Gamma$
having curvature $\kappa(\zeta)\ge \kappa$ with a certain
constant $\kappa>0$ and for all points $\zeta\in\Gamma$. Then
$M(K)\ge (\kappa/2) n$.
\end{theorem}

\begin{proof} As in \cite{R}, our proof hinges upon geometry in
a large extent. For this smooth case we use the following result,
which is well-known as Blaschke's Rolling Ball Theorem, cf.
\cite[p. 116]{Bla}.

\begin{lemma}[Blaschke]\label{l:Blaschke} Assume that the convex
domain $K$ has $C^2$ boundary $\Gamma=\partial K$ and that there
exists a positive constant $\kappa>0$ such that the curvature
$\kappa(\zeta)\ge \kappa$ at all boundary points
$\zeta\in\Gamma$. Then to each boundary points $\zeta\in\Gamma$
there exists a disk $D_R$ of radius $R=1/\kappa$, such that
$\zeta\in\partial D_R$, and $K\subset D_R$.
\end{lemma}

That is, if the curvature of the boundary curve of a twice
differentiable convex body exceeds $1/R$, then the convex body is
$R$-circular. From this an application of Theorem
\ref{th:circular} yields the assertion.
\end{proof}

So now it is worthy to calculate the curvature of $\partial E_b$.

\begin{lemma}\label{l:ellipsecurvature}
Let $E_b$ be the ellipse with major axes $[-1,1]$ and minor axes $[-ib,ib]$.
Consider its boundary curve $\Gamma_b$. Then at any point of the curve the
curvature is between $b$ and $1/b^2$.
\end{lemma}

\begin{proof} Now we depart from arc length parameterization and use for
$\Gamma_b:=\partial E_b$ the parameterization
$\gamma(\varphi):=(\cos(\varphi),b\sin(\varphi))$. Then we have
$$
\kappa(\gamma(\varphi))=
\frac{|\dot{\gamma}(\varphi)\times\ddot{\gamma}(\varphi)|}
{|\dot{\gamma}(\varphi)|^3}~,
$$
that is,

\begin{align}\label{ellipsecurvature}
\kappa(\gamma(\varphi)) & = \frac{\left|(-\sin\varphi,b\cos\varphi)
\times(-\cos\varphi,-b\sin\varphi)\right|} {|(-\sin\varphi,b\cos\varphi)|^3}
\notag \\
& = \frac{b\sin^2\varphi+b\cos^2\varphi}
{(\sin^2\varphi+b^2\cos^2\varphi)^{3/2}} \notag \\
& = \frac {b}{(\sin^2\varphi+b^2\cos^2\varphi)^{3/2}}~. \notag
\end{align}

Clearly, the denominator falls between
$(b^2\sin^2\varphi+b^2\cos^2\varphi)^{3/2}=b^3$ and
$(\sin^2\varphi+\cos^2\varphi)^{3/2}=1$, and these bounds are
attained, hence $\kappa(\gamma(\varphi))\in [b,1/b^2]$ whenever
$b\le 1$.
\end{proof}

\begin{proof}[Proof of Theorem \ref{th:ellipse}]
The curvature of $\Gamma_b$ at any of its points is at least $b$
according to Lemma \ref{l:ellipsecurvature}. Hence $M(E_b)\ge
(b/2) n$ in view of Theorem \ref{th:twicesmooth}, and Theorem
\ref{th:ellipse} follows.
\end{proof}

However, not only smooth convex domains can be proved to be
circular. Eg. it is easy to see that if a domain is the
intersection of finitely many $R$-circular domains, then it is
also $R$-circular. The next generalization is not that simple, but
is still true.

\begin{lemma}[Stranzen]\label{l:roughcurvature} Let the convex
domain $K$ have boundary $\Gamma=\partial K$ with angle function
$\alpha_{\pm}$ and let $\kappa>0$ be a fixed constant. Assume that
$\alpha_{\pm}$ satisfies the curvature condition
$\kappa(s)=\alpha'(s)\geq \kappa$ almost everywhere. Then to each
boundary point $\zeta\in\Gamma$ there exists a disk $D_R$ of
radius $R=1/\kappa$, such that $\zeta\in\partial D_R$, and
$K\subset D_R$. That is, $K$ is $R=1/\kappa$-circular.
\end{lemma}

\begin{proof} This result is essentially the far-reaching, relatively
recent generalization of Blaschke's Rolling Ball Theorem by
Stranzen. A reference for it is Lemma 9.11 on p. 83 of \cite{BS}.
Note that the proof of this lemma starts with establishing
Condition (i) on p. 83 of \cite{BS}, which is equivalent to the
subdifferential condition \eqref{subdiffcond}. We could as well
choose any of the equivalent formulations in
\eqref{aedifferentiability}-\eqref{subdiffcond}. The only slight
alteration from the formulation, suppressed in the above
quotations, is that Stranzen's version assumes $\kappa(t)\geq
\kappa$ wherever the curvature $\kappa(t)=\alpha'(t)$ exists (so
almost everywhere for sure), while above we stated the same thing
for almost everywhere, but not necessarily at every points of
existence. This can be overcome by reference to the
subdifferential version, too. Also, there is an even more recent
proof, which provides this version directly, see \cite{RB}.
\end{proof}

\begin{theorem}\label{th:aecurvature} Assume that the convex domain
$K$ has boundary $\Gamma=\partial K$ and that the a.e. existing
curvature of $\Gamma$ exceeds $\kappa$ almost everywhere, or,
equivalently, assume the subdifferential condition
\eqref{subdiffcond} (or any of the equivalent formulations in
\eqref{aedifferentiability}-\eqref{subdiffcond}) with
$\lambda=\kappa$. Then for all $p\in\PK$ we have
\begin{equation}\label{subdiffresult}
\|p'\| \ge \frac{\kappa}{2} n \|p\|~.
\end{equation}
\end{theorem}
\begin{proof} The proof follows from a combination of
Theorem \ref{th:circular} and Lemma \ref{l:roughcurvature}.
\end{proof}


Let us illustrate the strengths and weaknesses of the above
results on the folowing instructive examples, suggested to us by
J. Szabados (personal communication). Consider for any $1 < p<
\infty$ the $\ell_p$ unit ball
\begin{equation}\label{lpunitball}
B^p:=\{(x,y)\,:\, |x|^p+|y|^p\leq 1 \},\qquad\qquad
\Gamma^p:=\partial B^p = \{(x,y)\,:\, |x|^p+|y|^p= 1 \}.
\end{equation}
Also, let us consider for any parameter $0<b\leq 1$ the affine
image ("$\ell_p$-ellipse")
\begin{equation}\label{lpellipse}
B_b^p:=\{(x,y)\,:\, |x|^p+|y/b|^p\leq 1 \},\qquad\quad
\Gamma_b^p:=\partial B_b^p = \{(x,y)\,:\, |x|^p+|y/b|^p= 1 \}.
\end{equation}

By symmetry, it suffices to analyze the boundary curve
$\Gamma:=\Gamma_b^p$ in the positive quadrant. Here it has a
parametrization $\Gamma(x):=(x,y(x))$, where $y(x)=b\left(1-x^p
\right)^{1/p}$. As above, the curvature of the general point of
the arc in the positive quadrant can be calculated and we get
\begin{equation}\label{gammacurvature}
\kappa(x)=\frac{(p-1)bx^{p-2}(1-x^p)^{1/p-2}}
{\left(1+b^2x^{2p-2}(1-x^p)^{2/p-2}\right)^{3/2}}
\end{equation}

For $p>2$, the curvature is continuous, but it does not stay off
$0$: e.g. at the upper point $x=0$ it vanishes. Therefore, neither
Theorem \ref{th:twicesmooth} nor Theorem \ref{th:aecurvature} can
provide any bound, while Theorem \ref{th:cn} provides an estimate,
even if with a small constant: here $d(B)=2$, $w(B)=2b$, and we
get $M(B)\geq 0.00015 b n$.

When $p=2$, we get back the disk and the ellipses: the curvature
is minimal at $\pm ib$, and its value is $b$ there, hence
$M(B)\geq (b/2)n$, as already seen in Theorem \ref{th:ellipse}. On
the other hand Theorem \ref{th:cn} yields only $M(B)\geq
0.00015bn$ also here.

For $1<p<2$ the situation changes: the curvature becomes infinite
at the "vertices" at $\pm ib$ and $\pm 1$, and the curvature has a
positive minimum over the curve $\Gamma$. When $b=1$, it is
possible to explicitly calculate it, since the role of $x$ and $y$
is symmetric in this case and it is natural to conjecture that
minimal curvature occurs at $y=x$; using geometric-arithmetic mean
and also the inequality between power means (i.e.
Cauchy-Schwartz), it is not hard to compute $\min \kappa(x,y) =
(p-1)2^{1/p-1/2}$, (which is the value attained at $y=x$). Hence
Theorem \ref{th:aecurvature} (but not Theorem
\ref{th:twicesmooth}, which assumes $C^2$-smoothness, violated
here at the vertices!) provides $M(B^p)\geq (p-1)2^{1/p-3/2} n$,
while Theorem \ref{th:cn} provides, in view of
$w(B^p)=2^{3/2-1/p}$, something like $M(B^p)\geq 0.0003
\,2^{-1/2-1/p}n\geq 0.0001 n$, which is much smaller until $p$
comes down very close to $1$.

For general $0<b<1$ we obviously have $d(B)=2$, $(\sqrt{2}b <)
2b/\sqrt{1+b^2}<w(B)<2b$, and Theorem \ref{th:cn} yields $M(B)\geq
0.0001bn$ independently of the value of $p$.

Now $\min \kappa$ can be estimated within a constant factor
(actually, when $b\to 0$, even asymptotically precisely) the
following way. On the one hand, taking $x_0:=2^{-1/p}$ leads to
$\kappa(x_0)=(p-1)b 2^{1+1/p}/(1+b^2)^{3/2}<b(p-1)2^{1+1/p}$. On
the other hand denoting $\xi:=x^p$ and $\beta:=2/p-1\in(0,1)$,
from \eqref{gammacurvature} we get
$$
\frac{(p-1)b}{\kappa(x)}=\left[\xi(1-\xi)\right]^{\beta}
\left[\xi^{1-\beta}+b^2(1-\xi)^{1-\beta}\right]^{3/2} \leq
2^{-2\beta}\left[(\xi+(1-\xi))^{1-\beta}
(1+(b^2)^{1/\beta})^{\beta}\right]^{3/2},
$$
with an application of geometric-arithmetic mean inequality in the
first and H\"older inequality in the second factor. Note that when
$b\to 0$, this is asymptotically equivalent to $\kappa(x_0)\sim
b(p-1)2^{1+1/p}$. In general we can just use $b<1$ and get
$$
{\kappa(x)}\geq {(p-1)b} 2^{2\beta} \left[
1+b^{2/\beta}\right]^{-3\beta/2}\geq {(p-1)b} 2^{\beta/2}=(p-1)b
2^{1/p-1/2},
$$
within a factor $2^{3/2}$ of the upper estimate for $\min \kappa$.

Therefore, inserting this into Theorem \ref{th:aecurvature} as
above, we derive $M(B^p_b)\geq (p-1)b 2^{1/p-3/2} n$.

In all, we see that Theorems \ref{th:twicesmooth} (essentially due
to Er\H od) and \ref{th:aecurvature} usually (but not always, c.f.
the case $p\approx 1$ above !)) give better constants, when they
apply. However, in cases the curvature is not bounded away from 0,
we can retreat to application to the fully general Theorem
\ref{th:cn}, which, even if with a small absolute constant factor,
but still gives a precise estimate even regarding dependence of
the constant on geometric features of the convex domain. This
latter phenomenon is not just an observation on some particular
examples, but is a general result, also proved in \cite{R}, valid
even for not necessarily convex domains.

\begin{othm}\label{oth:sharp} Let $K\subset\CC$ be any
compact, connected set with diameter $d$ and minimal width $w$.
Then for all $n>n_0:=n_0(K):= 2 (d/16w)^2 \log (d/16w)$ there
exists a polynomial $p\in\PK$ of degree exactly $n$ satisfying
\begin{equation}\label{sharpresult}
\Norm{p'} \leq ~ C'(K)~ n~ {\Norm{p}} \qquad \text{\rm with}
\qquad C'(K):= 600 ~\frac{w(K)}{d^2(K)}~.
\end{equation}
\end{othm}


\mysection{Further remarks and problems} \label{sec:questions}

In the case of the unit interval also Tur\'an type $L^p$ estimates
were studied, see \cite{Zhou} and the references therein. It would
be interesting to consider the analogous question for convex
domains on the plane. Note that already Tur\'an remarked, see the
footnote in \cite[p.141]{Tur}, that on $D$ an $L^p$ version holds,
too. Also note that for domains there are two possibilities for
taking integral norms, one being on the boundary curve and another
one of integrating with respect to area. It seems that the latter
is less appropriate and convenient here.

In the above we described a more or less satisfactory answer of
the problem of inverse Markov factors for convex domains. However,
Levenberg and Poletsky showed that starshaped domains already do
not admit similar inverse Markov factors. A question, posed by V.
Totik, is to determine exact order of the inverse Markov factor
for the "cross" $C:=[-1,1]\cup [-i,i]$; clearly, the point is not
in the answer for the cross itself, but in the description of the
inverse Markov factor for some more general classes of sets.

Another question, still open, stems from the Szeg\H o extension of
the Markov inequality, see \cite{Szeg}, to domains with sector
condition on their boundary. More precisely, at $z\in\partial K$
$K$ satisfies the \emph{outer sector condition} with $0<\beta<2$,
if there exists a small neighborhood of $z$ where some sector
$\{\zeta~:~ \arg(\zeta-z) \in (\theta,\beta\pi+\theta) \}$ is
disjoint from $K$. Szeg\H o proved, that if for a domain $K$,
bounded by finitely many smooth (analytic) Jordan arcs, the
supremum of $\beta$-values satisfying outer sector conditions at
some boundary point is $\alpha<2\pi$, then $\|P'\| \ll
n^{\alpha}\|P\|$ on $K$. Then Tur\'an writes: "Es ist sehr
wahrscheinlich, da{\ss} auch den Szeg\H oschen Bereichen $M(p)\geq
c n ^{1/\alpha}$...", that is, he finds it rather likely that the
natural converse inequality, suggested by the known cases of the
disk and the interval (and now also by any other convex domain)
holds also for general domains with outer sector conditions.

\mysection{Acknowledgement}

The author is indebted to J. Kincses, E. Makai and V. Totik for
useful discussions, in particular for calling his attention to
the references \cite{Bla} and \cite{BS}.

\noindent
\ \\
{\bf Bibliography}
\\

\noindent {\sc\small Alfr\' ed R\'enyi Institute of Mathematics, \\
Hungarian Academy of Sciences, \\
Re\'altanoda utca 13-15., \\
1054 Budapest, Hungary} \\
e-mail: {\tt revesz@renyi.hu}

and

\noindent {\sc\small Institute Henri Poincar\'e, \\
11 rue Pierre et Marie Curie \\
75005 Paris, France} \\
e-mail: {\tt Szilard.Revesz@ihp.jussieu.fr}


\begin{thebibliography}{AAA}
\vspace{-1.6cm}

\bibitem{Bla} W. Blaschke, {\em Kreis und Kugel}, Zweite
Auflage, Walter de Gruyter AG, Berlin, 1956.

\bibitem{BS} J. N. Brooks, J. B. Stranzen, Blaschke's rolling ball theorem
in $\RR^n$, {\em Mem. Amer. Math. Soc.} {\bf 80, \# 405}, American
Mathematical Society, 1989.

\bibitem{E2} T. Erd\'elyi, Inequalities for exponential sums via
interpolation and Tur\'an type reverse Markov inequalities, {\em
manuscript}, 2005.
\texttt{www.math.tamu.edu/$\sim$tamas.erdelyi/papers-online/SHARMA\_sub.pdf}.

\bibitem{Er} J. Er\H od, Bizonyos polinomok maximum\'anak
als\'o korl\'atj\'ar\'ol, {\em Mat. Fiz. Lapok}\ {\bf 46} (1939),
58-82 (in Hungarian).

\bibitem{Faber} G. Faber, \"Uber Tschebyscheffsche Polynome, {\em
J. Reine Angew. Math.} {\bf 150} (1919), 79--106.

\bibitem{LP} N. Levenberg, E. Poletsky, Reverse Markov
inequalities, {\em Ann. Acad. Fenn.} {\bf 27} (2002), 173-182.


\bibitem{MMR} G. V. Milovanovi\'c, D. S. Mitrinovi\'c, Th. M. Rassias,
{\em Topics in Polynomials: Extremal Problems, Inequalities,
Zeros}, World Scientific, Singapore, 1994.

\bibitem{MR} G. V. Milovanovi\'c, Th. M. Rassias, New developments
on Tur\'an's extremal problems for polynomials, in:
\emph{Approximation Theory: In Memoriam A. K. Varma}, Marcel
Decker Inc., New York, 1998, pp. 433-447.

\bibitem{Rans} T. Ransford, {\em Potential Theory in
the Complex Plane}, London Mathematical Society Student Texts {\bf
28}, Cambridge University Press, 1994.

\bibitem{R} Sz. Gy. R\'ev\'esz, Tur\'an-type converse Markov inequalities
for convex domains on the plane, \emph{J. Approx. Theory}, to
appear. \comment{{\em Preprint of the Alfr\'ed R\'enyi Institute
of Mathematics}, {\bf \# 4/2005}.
\texttt{http://arxiv.org/abs/math.CA/0504416}}

\bibitem{RB} Sz. Gy. R\'ev\'esz, A discrete extension of the Blaschke Rolling Ball
Theorem, \emph{manuscript} \comment{{\em Preprint of the Alfr\'ed
R\'enyi Institute of Mathematics}, {\bf \# XXX/2006}.
\texttt{http://arxiv.org/abs/math.YYYYY/00000000}}.

\bibitem{Riesz} M. Riesz, Eine trigonometrische
Interpolationsformel und einige Ungleichungen f\"ur Polynome, {\em
Jahrsber. der deutsher Math. Vereinigung}, {\bf 23}, (1914),
354--368.

\bibitem{Sew} W. E. Sewell, On the polynomial derivative constant
for an ellipse, {\em Amer. Math. Monthly}, {\bf 44} (1937),
577-578.

\bibitem{Szeg} G. Szeg\H o, \"Uber einen Satz von A. Markoff, {\em
Math. Zeitschrift} {\bf 23} (1923), 45--61.

\bibitem{Tur} P. Tur\'an, \"Uber die Ableitung von Polynomen,
{\em Comp. Math.}\ {\bf 7} (1939), 89-95.

\bibitem{Zhou} S. P. Zhou, Some remarks on Tur\'an's inequality III:
the completion, {\em Anal. Math.} {\bf 21} (1995), 313-318.

\end{thebibliography}
\end{document}